\documentclass[epsfig,latexsym,amsfonts,twoside]{article}
\usepackage{amssymb,amsmath,amscd}

\def\part#1{\frac{\partial\phantom{#1}}{\partial#1}}
\newtheorem{thm}{Theorem}

\newtheorem{lem}[thm]{Lemma}

\newtheorem{cnj}[thm]{Conjecture}
\newenvironment{prf}{\begin{trivlist}\item[]{\bf Proof} }%
{\hfill $\Box$ \end{trivlist}}
{\end{trivlist}}
\newenvironment{rmk}{\begin{trivlist}\item[]{\bf Remark} }%
{\end{trivlist}}
{\end{trivlist}}

\def\Z{\ifmmode{{\mathbb Z}}\else{${\mathbb Z}$}\fi}
\def\Q{\ifmmode{{\mathbb Q}}\else{${\mathbb Q}$}\fi}
\def\C{\ifmmode{{\mathbb C}}\else{${\mathbb C}$}\fi} 
\def\P{\ifmmode{{\mathbb P}}\else{${\mathbb P}$}\fi} 

\def\H{\ifmmode{{\mathrm H}}\else{${\mathrm H}$}\fi} 

\def\B{\ifmmode{{\cal B}}\else{${\cal B}$}\fi} 
\def\E{\ifmmode{{\cal E}}\else{${\cal E}$}\fi} 
\def\F{\ifmmode{{\cal F}}\else{${\cal F}$}\fi} 
\def\K{\ifmmode{{\cal K}}\else{${\cal K}$}\fi} 
\def\L{\ifmmode{{\cal L}}\else{${\cal L}$}\fi} 
\def\M{\ifmmode{{\cal M}}\else{${\cal M}$}\fi} 
\def\N{\ifmmode{{\cal N}}\else{${\cal N}$}\fi} 
\def\O{\ifmmode{{\cal O}}\else{${\cal O}$}\fi} 
\def\U{\ifmmode{{\cal U}}\else{${\cal U}$}\fi}
\def\X{\ifmmode{{\cal X}}\else{${\cal X}$}\fi} 

\def\Br{\ifmmode{{\mathrm{Br}}}\else{${\mathrm{Br}}$}\fi} 
\def\OG{\ifmmode{\widetilde{\cal M}_4}\else{$\widetilde{\cal M}_4$}\fi} 
\def\D{\ifmmode{{\cal{D}}_{\mathrm{coh}}^b}\else{${{\cal{D}}_{\mathrm{coh}}^b}$}\fi}
\def\Shah{\ifmmode{\amalg\hspace*{-3.5pt}\amalg}\else{$\amalg\hspace*{-3.5pt}\amalg$}\fi}

\begin{document}

\title{Lagrangian fibrations on Hilbert schemes of points on K3
  surfaces} 
\author{Justin Sawon}
\date{October, 2005}
\maketitle

\begin{abstract}
Let $\mathrm{Hilb}^gS$ be the Hilbert scheme of $g$ points on a K3
surface $S$. Suppose that $\mathrm{Pic}S\cong\Z C$ where $C$ is a
smooth curve with $C^2=2(g-1)n^2$. We prove that $\mathrm{Hilb}^gS$
is a Lagrangian fibration.
\end{abstract}

\section{Introduction}

In this article we address the question: when does a holomorphic
symplectic manifold admit a fibration by Lagrangian tori? 

Suppose we have a (proper, with connected fibres) surjective morphism
$\pi:X\rightarrow Y$ between smooth projective varieties. Let $H_X$
and $H_Y$ denote ample line bundles on $X$ and $Y$
respectively. Consider the nef line bundle $L:=\pi^*H_Y$ on $X$; $L$
has numerical Kodaira dimension $\mathrm{dim}Y$, meaning that this is
the largest integer $j$ such that the cup-product of $j$ copies of
$L$ (thought of as a class in $\H^2(X,\Z)$) is non-trivial in
$\H^{2j}(X,\Z)$. Another way of stating this is that the polynomial
$$\int_X(H_X+tL)^{\mathrm{dim}X}$$
in $t$ has degree $\mathrm{dim}Y$.

The above observations are elementary, but now suppose that $X$ is an
irreducible holomorphic symplectic manifold of dimension $2g$,
i.e. $X$ is compact, K{\"a}hler, and $\H^0(X,\Omega^2)$ is generated
by a closed two-form $\sigma$ which is non-degenerate (meaning it
induces an isomorphism $T\cong\Omega^1$, or equivalently,
$\sigma^{\wedge g}$ trivializes the canonical bundle $K_X$). It is a
deep fact about holomorphic symplectic manifolds that the degree $2g$
form $\int_X\phi^{2g}$ on $\H^2(X,\Z)$ is the $g$th power of a
quadratic form, up to a constant. The quadratic form $q_X$ is known as
the Beauville-Bogomolov form~\cite{beauville83}, and
$$\int_X\phi^{2g}=c_Xq_X(\phi)^g$$
is the Fujiki relation. In particular
$$\int_X(H_X+tL)^{2g}=c_Xq_X(H_X+tL)^g$$
must be of degree $0$, $g$, or $2g$. This is essentially how
Matsushita~\cite{matsushita99,matsushita00i} proved that the base and
generic fibre of a non-trivial surjection from $X$ must be of
dimension $g$. By including also the term $s\sigma$ in the Fujiki
relation, he showed that the fibres are (holomorphic) Lagrangian,
and in particular the generic fibre must be a complex torus.

When we have such a Lagrangian fibration on $X$, the line bundle $L$
is isotropic with respect to $q_X$. Conversely
\begin{cnj}[Huybrechts~\cite{ghj02}, Sawon~\cite{sawon03}]
If there exists a non-trivial line bundle $L$ on $X$ which satisfies
$q_X(L)=0$, then $X$ is a Lagrangian fibration, at least up to
birational equivalence.
\end{cnj}
Associated to $L$ there are rational maps
$$\phi_{L^{\otimes N}}:X\dashrightarrow \P\H^0(X,L^{\otimes
  N})^{\vee}$$
and part of the problem is to show that $L^{\otimes N}$ is base-point
free, so that the map is a morphism. The Kodaira dimension
$$\lim_{N\rightarrow\infty}\mathrm{dim}\mathrm{Im}\phi_{L^{\otimes
    N}}$$
of a line bundle is never greater than its numerical Kodaira
dimension, and the latter is equal to $g$ when $L$ is isotropic by the
remarks above; thus the other part of the problem is to show that we
actually have equality here. This amounts to showing that $L$ has
sufficiently many sections and hence is a special case of the
Abundance Conjecture. However, in this article we will approach the
conjecture in a less direct way.

An example due to Beauville~\cite{beauville83} of an irreducible
holomorphic symplectic manifold is the Hilbert scheme of $g$ points on
a K3 surface $S$, denoted $\mathrm{Hilb}^gS$. In this case there is an
isometric isomorphism of weight-two Hodge structures
$$(\H^2(\mathrm{Hilb}^gS,\Z),q_X)=\H^2(S,\Z)\oplus\Z E$$
where $\H^2(S,\Z)$ has the usual intersection pairing, $E$ has square
$-2(g-1)$, and the direct sum is an orthogonal decomposition. Suppose
$L=(C,nE)$ is isotropic; then
$$C^2=2(g-1)n^2.$$
If $n=0$, then $C^2=0$, which implies that the K3 surface $S$ is
elliptic. The elliptic fibration then induces a Lagrangian fibration
on the Hilbert scheme, i.e.\
$$\mathrm{Hilb}^gS\rightarrow\mathrm{Hilb}^g\P^1=\mathrm{Sym}^g\P^1\cong\P^g.$$  
For $n\geq 1$, assume that $\mathrm{Pic}S\cong\Z C$. In particular, we
assume $C$ is primitive, so it must be reduced and irreducible. If
$n=1$, then $C^2=2g-2$, and $C$ will represent a smooth genus $g$
curve. One can show that $\mathrm{Hilb}^gS$ is birational to the
compactified Jacobian of the linear system $|C|\cong\P^g$. More
precisely, denote the family of curves by ${\cal C}/\P^g$; then
$\mathrm{Hilb}^gS$ is birational to $\overline{\mathrm{Pic}}^g({\cal
  C}/\P^g)$. This example was studied by
Beauville~\cite{beauville99}. The case $g=2$ and $n=2$ also leads to a
fibration (see Hassett and Tschinkel~\cite{ht99} or Fu~\cite{fu03}),
which we describe in Section 2. In fact, using a deformation argument
Hassett and Tschinkel proved that $\mathrm{Hilb}^2S$ admits a
fibration when $S$ is a generic K3 surface containing a smooth curve
$C$ with $C^2=2n^2$ (the case $g=2$ and $n$ arbitrary). More recently
Iliev and Ranestad~\cite{ir05} constructed a fibration on
$\mathrm{Hilb}^3S$ when $C$ is a genus 9 curve (the case $g=3$ and
$n=2$).

In this article we prove 
\begin{thm}
\label{main_thm}
Let $S$ be a K3 surface with $\mathrm{Pic}S\cong\Z C$ where $C$ is a
smooth curve with $C^2=2(g-1)n^2$ ($g\geq 2$ and $n\geq 2$). Then
$\mathrm{Hilb}^gS$ is a Lagrangian fibration.
\end{thm}
In general we expect that the fibration should be the compactified
Jacobian of a family of curves, but the K3 surface $S$ does not
contain curves of genus $g$. However, by studying carefully the case
$g=2$ and $n=2$ we see that there is a Mukai dual K3 surface
$S^{\prime}$ which does contain curves of genus $g$. We then use a
twisted Fourier-Mukai transform to prove that $\mathrm{Hilb}^gS$ is
isomorphic to the compactified Jacobian of the curves on $S^{\prime}$
(more precisely, to a torsor over the compactified Jacobian); in
particular, it is a Lagrangian fibration. A similar result has been
obtained independently and simultaneously by
Markushevich~\cite{markushevich05}, also by using Fourier-Mukai
transforms.

Unfortunately the proofs of some of our lemmas fail when $C$ is not
primitive, though we expect the theorem should still be true. In the
final section we outline how our method ought to generalize from
Hilbert schemes to arbitrary moduli spaces of sheaves on K3 surfaces.

The author would like to thank Claudio Bartocci, Tom Bridgeland,
Daniel Huybrechts, Manfred Lehn, and Paolo Stellari for useful
conversations, and Ivan-penguin for help in going from big to
ample. Thanks also to Kota Yoshioka and Dimitri Markushevich for
pointing out an error in an earlier version of this paper, and
especially to Kota Yoshioka for help in correcting it. The author is
grateful for the hospitality of the Johannes-Gutenberg
Universit{\"a}t, Mainz, where this article was written. The author is
supported by NSF grant number 0305865.

\section{A motivating example}

We begin by describing a Lagrangian fibration on $\mathrm{Hilb}^2S$
for a particular K3 surface $S$. This example was described by Hassett
and Tschinkel~\cite{ht99} and by Fu~\cite{fu03}, though it appears to
have been known to experts for some time before those articles. The
fibration arises somewhat unexpectedly, but leads to a very rich
geometry.

\subsection{Classical description}

Let $S$ be a K3 surface given by the intersection of three generic
smooth quadrics $Q_1$, $Q_2$, and $Q_3$ in $\P^5$. Then $S$ contains a
genus five curve $C$ with $C^2=8$, and by the genericity assumption
$\mathrm{Pic}S$ is generated by $C$ (we will tend to use the same
notation for the curve $C$, its divisor class, and its class in
$\mathrm{H}^2(S,\Z)$; the meaning will always be obvious from the
context). Let $\xi$ be a zero-dimensional
length two subscheme of $S$. Then $\xi$ is either a pair of distinct
points or a point with a tangent direction, and in each case we obtain
a line $\ell$ in $\P^5$. Each quadric in the net $(\P^2)^{\vee}$
generated by $Q_1$, $Q_2$, and $Q_3$ contains $S$, and therefore
intersects $\ell$ in $\xi$. If a quadric intersects $\ell$ in one
additional point, then it must contain $\ell$. This is a single linear
condition, and therefore quadrics which contain $\ell$ form a pencil
$p\cong\P^1$ in $(\P^2)^{\vee}$. Dually, $\xi$ determines a point $P$
in $\P^2$ and we obtain a morphism
$$\mathrm{Hilb}^2S\rightarrow\P^2.$$
By general results of Matsushita~\cite{matsushita99,matsushita00i} on
the structure of fibrations on irreducible holomorphic symplectic
manifolds, we know that the generic fibre must be a Lagrangian abelian 
surface. Moreover we can give an explicit description.

Fix a point $P\in\P^2$ corresponding to a pencil $p$ of
quadrics. The base locus of this pencil is a three-fold $Y_P$ in
$\P^5$, and each line $\ell$ in $Y_P$ intersects $S$ in a length two
subscheme $\xi$. Therefore the fibre over $P$ is the space
parameterizing lines in $Y_P$, which is known to be an abelian surface
(see Chapter 6 of Griffiths and Harris~\cite{gh78}). In fact,
Griffiths and Harris give
a description of this surface as the Jacobian of a genus two curve: in
the pencil of quadrics $p$ there will be precisely six singular
quadrics (counted with multiplicity), and the genus two curve $D_p$ is
the double cover of $p\cong\P^1$ branched over these six points. Note
that if some of these six points coincide then $D_p$ is singular; in
that case the fibre is the compactified Jacobian of
D'Souza~\cite{dsouza79}, which is a degeneration of a nearby smooth
fibre. Altogether we see that $\mathrm{Hilb}^2S$ is locally isomorphic
as a fibration to the relative compactified Jacobian of the family of
curves
$$\{D_p|p\subset (\P^2)^{\vee}\mbox{ or equivalently }P\in\P^2\}.$$
We don't claim that the isomorphism is global: see our comments in the 
next subsection.

\subsection{Moduli space interpretation}

There is a result of Markushevich~\cite{markushevich95,markushevich96}
which roughly says the following: suppose we have a family of genus
two curves parameterized by a complex surface $B$ such that the
compactified Jacobian of the family is a holomorphic symplectic
four-fold; then the curves are the complete linear system of curves on
a K3 surface, and in particular $B$ is the projective plane. Curiously
in the above example $S$ does not contain any genus two
curves. However, we will describe another K3 surface $S^{\prime}$
which does, and indeed these are precisely the curves $D_p$ described
above. The K3 surface $S^{\prime}$ will be Mukai dual to $S$, and the
(twisted) Fourier-Mukai transform will induce the isomorphism between 
$\mathrm{Hilb}^2S$ and the relative Jacobian of the linear system of
curves $D_P$. More precisely, $\mathrm{Hilb}^2S$ will be isomorphic to
a moduli space of stable (torsion) twisted sheaves on $S^{\prime}$,
and the latter will only be locally isomorphic to the relative
Jacobian, as a fibration.

In fact $S^{\prime}$ is easy to describe. In the net $(\P^2)^{\vee}$
of quadrics there is a sextic parameterizing singular quadrics, and
$S^{\prime}$ is the double cover of $(\P^2)^{\vee}$ branched over this
sextic. Clearly it contains the curves $D_p$, which are double covers
of the pencils in $p\subset (\P^2)^{\vee}$ branched over the same
points. 

We want to understand $S^{\prime}$ as a moduli space of sheaves on
$S$. This is Example 0.9 in Mukai~\cite{mukai84} (see also Example 2.2
in~\cite{mukai88} for more details). Each quadric in the
net $(\P^2)^{\vee}$ can be regarded as an embedding of the complex
Grassmannian $\mathrm{Gr}(2,4)$ of planes in $\C^4$ into $\P^5$
(strictly speaking it is a birational map of $\mathrm{Gr}(2,4)$ onto
its image when the quadric is singular). Over the Grassmannian
there are two natural rank two bundles: the universal bundle of planes
$E$, which sits inside the trivial rank four bundle, and the universal
quotient bundle $F$. Thinking of $E^{\vee}$ and $F$ as bundles over
the quadric, and restricting to the K3 surface $S$, gives two stable
bundles over $S$ with $c_1=C$ and $c_2=4$. In other words, they have
Mukai vector
$$w:=(2,C,2).$$
In case the quadric is singular, $E^{\vee}|_S$ and $F|_S$ become
isomorphic. Thus $S^{\prime}$, the double cover of $(\P^2)^{\vee}$
branched over the sextic of singular quadrics, parameterizes a family
of stable bundles on $S$ of Mukai vector $w$. Since $w$ is isotropic,
we know that the moduli space $M_w(S)$ has dimension two, and thus
$$S^{\prime}=M_w(S).$$

The Hilbert scheme $\mathrm{Hilb}^2S$ can be regarded as the moduli
space $M_{(1,0,-1)}(S)$ of stable sheaves on $S$ with Mukai vector
$v:=(1,0,-1)$. The degree $d$ relative compactified Picard scheme
$\overline{\mathrm{Pic}}^d$ of the family of genus two curves on
$S^{\prime}$ can be regarded as the moduli space
$M_{(0,D,d-1)}(S^{\prime})$ of stable sheaves on $S^{\prime}$ with
Mukai vector $v^{\prime}:=(0,D,d-1)$ (these are torsion sheaves
supported on a curve in the linear system $|D|$). We'd like to obtain 
an isomorphism
$$M_{(1,0,-1)}(S)\cong M_{(0,D,d-1)}(S^{\prime})$$
for some $d\in\Z$ from a Fourier-Mukai transform between $S$ and
$S^{\prime}=M_w(S)$. However, looking at $w=(2,C,2)$, we see that the
greatest common divisor of $2$, $\mathrm{deg}C=8$, and $2$ is not
one. Therefore $S^{\prime}$ is a non-fine moduli space of sheaves on
$S$, and a universal sheaf does not exist: the obstruction is a gerbe
$\alpha\in\mathrm{H}^2(S^{\prime},\O^*)$ (we will see later that
$\alpha$ really is non-trivial). However, in this situation there
exists a twisted universal sheaf, and this induces a twisted
Fourier-Mukai transform (see C{\u a}ld{\u a}raru~\cite{caldararu00})
$$\Phi_{S^{\prime}\rightarrow
  S}:\D(S^{\prime},\alpha^{-1})\stackrel{\cong}{\longrightarrow}\D(S).$$
The twisted Fourier-Mukai transform will induce an isomorphism
$$M_{(1,0,-1)}(S)\cong M_{(0,D,d-1)}(S^{\prime},\alpha^{-1})$$
between $\mathrm{Hilb}^2S$ and the moduli space of stable twisted
sheaves on $S^{\prime}$ with twisted Mukai vector $(0,D,d-1)$ (moduli
spaces of stable twisted sheaves were constructed by
Yoshioka~\cite{yoshioka04}). Note that because the latter are
supported on curves, and any gerbe restricted to a curve becomes
trivializable, any twisted sheaf in
$M_{(0,D,d-1)}(S^{\prime},\alpha^{-1})$ can be identified with an
untwisted sheaf. Because the identification is not canonical
$M_{(0,D,d-1)}(S^{\prime},\alpha^{-1})$ will not be globally
isomorphic to the untwisted moduli space
$M_{(0,D,d-1)}(S^{\prime})$. However, they are both torsors over the
relative compactified Jacobian of the family of curves, and in
particular $M_{(0,D,d-1)}(S^{\prime},\alpha^{-1})$ is a fibration over 
$|D|\cong\P^2$, the map given by taking support.

Since we have outlined the basic ideas, we will pass now to the
general situation, rather than working through this example in
detail.

\section{A Mukai dual K3 surface}

Let $S$ be a K3 surface and let $\mathrm{Hilb}^gS=M_{(1,0,1-g)}(S)$ be
the Hilbert scheme of $g\geq 2$ points on $S$, or equivalently the
moduli space of stable sheaves with Mukai vector $v:=(1,0,1-g)$. The
weight-two Hodge structure of moduli spaces of sheaves on a K3 surface
was described by O'Grady~\cite{ogrady97}. We have
$$\H^2(M_{(1,0,1-g)}(S),\Z)=(1,0,1-g)^{\perp}=\H^2(S,\Z)\oplus
\Z E$$
where $E^2=-2(g-1)$. Suppose there is a divisor on $\mathrm{Hilb}^nS$
with square zero with respect to the Beauville-Bogomolov quadratic
form. Write this divisor as $(C,nE)$ where
$$C^2=2(g-1)n^2.$$
We can take $n$ to be positive; in fact let $n\geq 2$, since the
existence of a Lagrangian fibration when $n=1$ is well-known (see the
introduction). We assume that our K3 surface $S$ has Picard number one
and that $\mathrm{Pic}S\cong\Z C$. Without loss of generality, suppose
that $C$ is ample (rather than $-C$); throughout, by $S$ we will
always mean $S$ with this given polarization. Since we've assumed $C$
is primitive, it must be reduced and irreducible.


Define a Mukai dual K3 surface $S^{\prime}$ as the moduli space
$M_w(S)$ of stable sheaves on $S$ with Mukai vector
$w:=(n,C,(g-1)n)$. Observe that $w$ is primitive and isotropic so
$S^{\prime}$ really is a K3 surface (Theorem 1.4 in
Mukai~\cite{mukai87}; non-emptyness of moduli spaces of stable sheaves
on a K3 surface was proved in the general case by Yoshioka
in~\cite{yoshioka99}, though many special cases were known earlier,
including $M_w(S)$ presumably).

\begin{lem}
Every sheaf in $S^{\prime}$ is locally free.
\end{lem}

\begin{prf}
Let $\U_t$ be the sheaf on $S$ corresponding to $t\in
S^{\prime}$. Since $\U_t$ is stable, it is $\mu-$semi-stable. But
$\U_t$ has slope
$$\mu(\U_t)=\frac{C^2}{n}$$
and since $C$ generates $\mathrm{Pic}S$, a proper subsheaf can have
slope at most zero. Thus $\U_t$ is $\mu-$stable. Now Corollary 3.10 of 
Mukai~\cite{mukai87} implies that $\U_t$ is locally free. Explicitly, 
the double-dual $\U_t^{\vee\vee}$ will also be $\mu-$stable. If $\U_t$
is not locally free then there is an exact sequence
$$0\rightarrow\U_t\rightarrow\U_t^{\vee\vee}\rightarrow\O_Z\rightarrow
0,$$
where $Z$ is a zero-dimensional subscheme of $S$, and
$$v(\U_t^{\vee\vee})=(n,C,(g-1)n+\ell(Z))$$
has square $-2n\ell(Z)<-2$, contradicting the Bogomolov inequality
(for example, see Huybrechts and Lehn~\cite{hl97}).
\end{prf}

Since the greatest common divisor of $n$, $\mathrm{deg}C=2(g-1)n^2$,
and $(g-1)n$ is not one, we do not expect $S^{\prime}$ to be a fine
moduli space of sheaves on $S$. More precisely, there exists a gerbe
$\alpha\in\H^2(S^{\prime},\O^*)$ which is the obstruction to the
existence of a universal sheaf, and this gerbe is
$n-$torsion. (Later we will see that $\alpha$ really is
non-trivial.) Let $B\in\H^2(S^{\prime},\Q)$ be any rational B-field
lift of $\alpha$, in the sense of Huybrechts and
Stellari~\cite{hs04}. In other words, the $(0,2)-$part of $B$ is
mapped to $\alpha$ under the exponential map.

Denote by $q$ and $p$ the projections from $S\times S^{\prime}$ to
$S$ respectively $S^{\prime}$. There exists a $p^*\alpha-$twisted
universal sheaf $\U$ on $S\times S^{\prime}$, which is such that
$\U_t=\U|_{S\times t}$ is the sheaf on $S$ which the point $t\in 
S^{\prime}$ parameterizes. There is a twisted Fourier-Mukai transform
(see C{\u a}ld{\u a}raru~\cite{caldararu00})
$$\Phi:=\Phi^{\U}_{S^{\prime}\rightarrow
  S}:\D(S^{\prime},\alpha^{-1})\rightarrow\D(S)$$
with quasi-inverse
$$\Psi:=\Psi^{\U^{\vee}}_{S\rightarrow
  S^{\prime}}:\D(S)\rightarrow\D(S^{\prime},\alpha^{-1}).$$
(The composition of $\Phi$ and $\Psi$ is the identity shifted by 2.)
At the cohomological level we obtain an isomorphism of twisted
weight-two Hodge structures (see~\cite{hs04}) 
$$\Phi_*:\tilde{\H}(S^{\prime},-B,\Z)\rightarrow\tilde{\H}(S,\Z)$$
by using $v(\U)$, which also preserves the Mukai pairing. Its
inverse $\Psi_*$ is given by using $v(\U^{\vee})$.

\begin{lem}
\label{HodgeH2}
There is an isometric isomorphism
$$\H^2(S^{\prime},\Z)\cong w^{\perp}/w.$$
\end{lem}

\begin{prf}
This is just an adaptation of the proof in the untwisted case (see
Theorem 6.1.14 in Huybrechts and Lehn~\cite{hl97}). The isomorphism is
induced by
$$\tilde{\H}(S,\Q)\stackrel{\Psi_*}{\longrightarrow}\tilde{\H}(S^{\prime},-B,\Q)\stackrel{\exp(B)}{\longrightarrow}\tilde{\H}(S^{\prime},\Q).$$
Firstly $\Psi_*$ is given by 
$$v(\U^{\vee})\in\H^*(S\times S^{\prime},-p^*B)=\H^*(S)\otimes\H^*(S^{\prime},-B).$$
The $\H^*(S)\otimes\H^0(S^{\prime},-B)$ part is given by $w^{\vee}$
(almost by definition). Therefore if $c\in\tilde{H}(S,\Q)$ is
orthogonal to $w$ then $\Psi_*(c)$ will have $\H^0(S^{\prime},-B,\Q)$
component zero. Next observe that $w$ is mapped to $(0,0,1)$ in both
$\tilde{\H}(S^{\prime},-B,\Q)$ and $\tilde{\H}(S^{\prime},\Q)$ (again,
almost by definition), so we can quotient out to obtain the
isomorphism
$$(w^{\perp}\otimes\Q/\Q w)\cong\H^2(S^{\prime},\Q)$$
over $\Q$. Since $\Psi_*$ and the isometry 
$$\H^2(S^{\prime},-B,\Q)\stackrel{\exp(B)}{\longrightarrow}\H^2(S^{\prime},\Q)\oplus\H^4(S^{\prime},\Q)\stackrel{\mbox{pr}}{\longrightarrow}\H^2(S^{\prime},\Q)$$
induced by $\exp(B)$ are defined over $\Z$, we actually have an
isomorphism over $\Z$. (Note that $\exp(B)$ itself won't be defined
over $\Z$.)
\end{prf}

\begin{lem}
\label{genusg}
The K3 $S^{\prime}$ contains a curve $D$ of arithmetic genus $g$.
\end{lem}

\begin{prf}
Note that $(1,0,1-g)$ is orthogonal to $w=(n,C,(g-1)n)$ and
lies in the $(1,1)-$part of the Hodge decomposition, so by the
previous lemma it gives a divisor $D$ on $S^{\prime}$. In Section 5
of~\cite{hs04}, Huybrechts and Stellari proved that the Hodge isometry
coming from a twisted Fourier-Mukai transform associated to a moduli
space is orientation preserving. Using their results, one can deduce
that $D$ is in the positive cone of $S^{\prime}$. In fact, $D$ is
ample as $S^{\prime}$ must have the same Picard number as $S$, namely
one. Since $D^2=2g-2$, $D$ has arithmetic genus $g$.
\end{prf}

\begin{rmk}
We now argue that $\alpha$ is non-trivial. Suppose otherwise; first
note that the previous two lemmas would still apply. We would also
have a Fourier-Mukai transform between $S$ and $S^{\prime}$, which at
the cohomological level induces an isometry between generalized Picard
groups (see 3.6 of Orlov~\cite{orlov97}). However the quadratic forms
on 
$$\mathrm{Pic}(S,0)\cong U\oplus \Z C$$
and
$$\mathrm{Pic}(S^{\prime},0)\cong U\oplus \Z D$$
have different determinants
$$-2(g-1)n^2\neq -2(g-1).$$
(We will prove shortly that $\mathrm{Pic}S^{\prime}\cong\Z D$; in any
case, if $D=mD_0$ where $D_0$ is primitive then the determinant of
$\mathrm{Pic}(S^{\prime},0)$ would be $-2(g-1)/m^2$, and we would
still have a contradiction.) Moreover, one can show that the order of
$\alpha$ is precisely $n$.
\end{rmk}

\begin{rmk}
It follows from the proof of Lemma~\ref{genusg} that $\Psi_*$ must map
$(1,0,1-g)$ to $(0,D,k)$ for some $k\in\Z$. Since $\Psi$ takes the
structure sheaf ${\cal O}_s$ of a point $s\in S$ to
$\U^{\vee}|_{s\times S^{\prime}}$, which is a rank $n$ twisted sheaf,
$\Psi_*$ will also take $(0,0,1)$ to some
$(n,-E,l)\in\tilde{\H}(S^{\prime},-B,\Z)$. We summarize the action of
$\Psi_*$:
\begin{eqnarray*}
w=(n,C,(g-1)n) & \mapsto & (0,0,1) \\
(1,0,1-g) & \mapsto & (0,D,k) \\
(0,0,1) & \mapsto & (n,-E,l). \\
\end{eqnarray*}
\end{rmk}

\begin{rmk}
In fact, $\Psi_*$ induces a map which takes a divisor $H$ on $S$ to the
divisor $H^{\prime}$ on $S^{\prime}$ given by the middle component of
$$-\Psi_*(H+\frac{1}{n}\langle H,C\rangle (0,0,1))$$
($C$ appears in this formula as the middle component of
$w=(n,C,(g-1)n)$, and $n$ appears as the order of $\alpha$). Therefore
the polarization $C^{\prime}$ on $S^{\prime}$ which corresponds to the
polarization $C$ of $S$ is the middle component of
$$-\Psi_*(C+\frac{1}{n}\langle C,C\rangle (0,0,1)),$$
and a calculation shows that this is equal to $nD$. Henceforth we will
always assume that $S^{\prime}$ comes equipped with this polarization.
\end{rmk}

For $s\in S$, $\U_s=\U|_{s\times S^{\prime}}$ is an $\alpha-$twisted
locally free sheaf of rank $n$ on $S^{\prime}$. (The reader should
take care not to confuse $\U_s$ with the sheaf $\U_t=\U|_{S\times t}$
on $S$.) In order to calculate the degree of an $\alpha-$twisted sheaf
$\cal F$, we choose a locally free $\alpha-$twisted sheaf $G$ whose
rank is the order of $\alpha$, i.e.\ $n$. We then define
$$\mathrm{deg}_G({\cal F}):=\mathrm{deg}(G^{\vee}\otimes{\cal F})$$
where the degree of the (untwisted) sheaf $G^{\vee}\otimes{\cal F}$ is
calculated with respect to the polarization $nD$ of $S^{\prime}$ (see
Yoshioka~\cite{yoshioka04}). We will choose $G$ to be $\U_s$ for some
$s\in S$; of course $\mathrm{deg}_G$ is independent of which $s\in S$
we choose.

\begin{lem}
\label{Picard}
The Picard group of $S^{\prime}$ is $\mathrm{Pic}S^{\prime}\cong\Z D$.
\end{lem}

\begin{prf}
Suppose that $D$ is not primitive; thus $D=mD_0$ where $D_0$ is
primitive and $m>1$. Let $L$ be any line bundle on $D_0$. Then
$\iota_*L$ is an $\alpha-$twisted sheaf on $S^{\prime}$, where
$\iota:D_0\hookrightarrow S^{\prime}$ is the inclusion of $D_0$ in
$S^{\prime}$. The degree of $\iota_*L$ is
$$\mathrm{deg}_G(\iota_*L)=\mathrm{deg}(G^{\vee}\otimes{\iota_*L})=nD.nD_0=\frac{2(g-1)n^2}{m}.$$
However, The twisted Fourier-Mukai transform must preserve the smallest
positive degree of (twisted) sheaves. In other words, since the
smallest positive degree of a sheaf on $S$ is $C^2=2(g-1)n^2$, an
$\alpha-$twisted sheaf on $S^{\prime}$ cannot have degree
$\frac{2(g-1)n^2}{m}$. This contradiction proves the lemma.
\end{prf}

\begin{lem}
\label{smoothcurves}
The linear system $|D|$ on $S^{\prime}$ contains smooth
curves. Equivalently, we can assume $D$ itself is smooth.
\end{lem}

\begin{prf}
Since $\mathrm{Pic}S^{\prime}\cong\Z D$, every curve in $|D|$ is
reduced and irreducible. The following proof that $|D|$ is base-point
free is taken from M{\'e}rindol~\cite{merindol85}: the short exact
sequence 
$$0\rightarrow\O_{S^{\prime}}\rightarrow\O_{S^{\prime}}(D)\rightarrow\O_D(D)\cong\omega_D\rightarrow
0$$
leads to the long exact sequence
$$0\rightarrow\H^0(S^{\prime},\O_{S^{\prime}})\rightarrow\H^0(S^{\prime},\O_{S^{\prime}}(D))\rightarrow\H^0(D,\omega_D)\rightarrow\H^1(S^{\prime},\O_{S^{\prime}})=0.$$
Therefore if $\O_{S^{\prime}}(D)$ had base-points, so too would $\O_{S^{\prime}}(D)|_D\cong\omega_D$;
but $D$ has arithmetic genus $D^2/2+1\geq 2$, so $\omega_D$ is
base-point free. 

Thus $|D|$ is base-point free, and by Bertini's Theorem the generic
element is a smooth curve.
\end{prf}

\section{Birationality of moduli spaces}

Our goal now is to use the twisted Fourier-Mukai transform between 
$S$ and $S^{\prime}$ to induce an isomorphism of (twisted) moduli
spaces
$$\mathrm{Hilb}^gS=M_{(1,0,1-g)}(S)\cong
M_{(0,D,k)}(S^{\prime},\alpha^{-1}).$$
In fact we will redefine $\Psi$ to be the Fourier-Mukai transform 
$$\Psi^{\U}_{S\rightarrow S^{\prime}}:\D(S)\rightarrow
\D(S^{\prime},\alpha)$$
rather than $\Psi^{\U^{\vee}}_{S\rightarrow S^{\prime}}$. By changing
$\U^{\vee}$ to $\U$, we find that $v=(1,0,1-g)$ is taken to
$(0,-D,k)$. In this section we will show that the composition of
$\Psi$ with taking cohomology and the shift functor induces a
birational map
$$\mathrm{Hilb}^gS=M_{(1,0,1-g)}(S)\dashleftarrow\dashrightarrow
M_{(0,D,-k)}(S^{\prime},\alpha).$$
In the next section we will show that this map is regular, and
therefore an isomorphism.

\subsection{The basic argument}

Let $Z\in\mathrm{Hilb}^gS$ be a generic element of the Hilbert scheme
of $S$. Thus $Z$ consists of $g$ distinct points
$\{z_1,\ldots,z_g\}$. The cohomological twisted Fourier-Mukai
transform $\Psi_*$ takes $v=(1,0,1-g)\in\tilde{\H}(S,\Z)$ to
$(0,-D,k)\in\tilde{\H}(S^{\prime},B,\Z)$. Therefore the Fourier-Mukai
transform $\Psi$ takes the ideal sheaf ${\cal I}_Z$, which has Mukai
vector $v$, to an object in $\D(S^{\prime},\alpha)$ with twisted Mukai
vector $(0,-D,k)$. A priori the image is a complex of twisted sheaves,
though we will show that it has non-trivial cohomology only in degree
one; thus it is canonically isomorphic to a (twisted) sheaf in degree
one. From its twisted Mukai vector, we know that this sheaf must look
like $\iota_*L[1]$, where $\iota:D_0\hookrightarrow S^{\prime}$ is
the inclusion of a curve $D_0\in|D|$ into $S^{\prime}$, and $L$ is a
(twisted) line bundle of some degree on $D_0$. Composing with the
shift functor we get a typical element of
$M_{(0,D,-k)}(S^{\prime},\alpha)$.

Now $\Psi$ takes ${\cal I}_Z$ to a complex whose cohomology ${\bf
  R}^ip_*(\U\otimes q^*{\cal I}_Z)$ vanishes for $i<0$ and $i>2$. Our
aim is to show that in fact only the ${\bf R}^1$ term is
non-vanishing, and for this we consider the fibres over $t\in
S^{\prime}$. Start with the short exact sequence
$$0\rightarrow{\cal I}_Z\rightarrow\O_S\rightarrow\O_Z\rightarrow 0.$$
Tensoring with $\U_t=\U|_{S\times t}$ and taking the corresponding long
exact sequence gives
$$0\rightarrow\H^0(S,{\cal
  I}_Z\otimes\U_t)\longrightarrow\H^0(S,\U_t)\stackrel{\mathrm{ev}}{\longrightarrow}\H^0(Z,\U_t|_Z)\rightarrow\phantom{0}$$
$$\phantom{0}\rightarrow\H^1(S,{\cal
  I}_Z\otimes\U_t)\longrightarrow\H^1(S,\U_t)\longrightarrow\H^1(Z,\U_t|_Z)\rightarrow\phantom{0}$$
$$\phantom{0}\rightarrow\H^2(S,{\cal
  I}_Z\otimes\U_t)\longrightarrow\H^2(S,\U_t)\longrightarrow\H^2(Z,\U_t|_Z)\rightarrow 0.$$
In the third column $\H^0(Z,\U_t|_Z)$ has dimension $gn$, since $\U_t$
is a rank $n$ vector bundle and $Z$ consists of $g$ points; the higher
cohomology vanishes. We also show that the higher cohomology of $\U_t$
on $S$ vanishes.

\begin{lem}
The groups $\H^1(S,\U_t)$ and $\H^2(S,\U_t)$ vanish for all $t\in
S^{\prime}$.
\end{lem}

\begin{prf}
Recall that $\U_t$ is a $\mu-$stable locally free sheaf with Mukai
vector $w=(n,C,(g-1)n)$. Thus $\U_t^{\vee}$ is also $\mu-$stable, with
Mukai vector $(n,-C,(g-1)n)$ and hence negative slope. Therefore
$\U_t^{\vee}$ has no sections and
$$\H^2(S,\U_t)\cong\H^0(S,\U_t^{\vee})^{\vee}=0.$$

Suppose $\H^1(S,\U_t)$ does not vanish. Then neither does
$$\mathrm{Ext}^1(\O,\U_t^{\vee})\cong\H^1(S,\U_t^{\vee})\cong\H^1(S,\U_t)^{\vee}$$
and hence we have a non-trivial extension
$$0\rightarrow\U_t^{\vee}\rightarrow{\cal
  F}\stackrel{\pi}{\longrightarrow}\O\rightarrow 0.$$
The Mukai vector of $\cal F$ is
$$v({\cal F})=v(\U_t^{\vee})+v(\O)=(n+1,-C,(g-1)n+1).$$
If $\cal F$ were stable then
$$v({\cal F})^2=-2(ng+1)<-2$$
would contradict the Bogomolov inequality~\cite{hl97}. Therefore $\cal
F$ is 
unstable. A destabilizing subsheaf must have rank less than $n+1$ and
slope greater than $-\frac{1}{n+1}C^2$, and thus at least zero since
$C$ generates $\mathrm{Pic}S$. By standard arguments (see
Friedman~\cite{friedman98}) we can assume that $\O$ is a destabilizing
subsheaf. Explicitly, suppose ${\cal F}^{\prime}$ destabilizes $\cal
F$ and consider the composition
$${\cal F}^{\prime}\hookrightarrow{\cal F}\stackrel{\pi}{\longrightarrow}\O.$$
If this is not an isomorphism then its kernel destabilizes
$\U_t^{\vee}$; but if it is an isomorphism then the sequence defining
$\cal F$ splits, also a contradiction.
\end{prf}

We immediately conclude from the long exact sequence that
$\H^2(S,{\cal I}_Z\otimes \U_t)$ vanishes for all $t\in S$, and hence
$${\bf R}^2p_*(\U\otimes q^*{\cal I}_Z)=0$$
for all $Z\in\mathrm{Hilb}^gS$.

\begin{lem}
The (twisted) sheaf ${\bf R}^0p_*\U$ on $S^{\prime}$ is locally free
of rank $gn$.
\end{lem}

\begin{prf}
By the previous lemma, the dimension of the fibre
$$({\bf R}^0p_*\U)_t=\H^0(S,\U_t)$$
is independent of $t\in S^{\prime}$ and equal to
$$\chi(\U_t)=-\langle v(\U_t),(1,0,1)\rangle=gn.$$
\end{prf}

Since ${\bf R}^0p_*(\U\otimes q^*{\cal I}_Z)$ is a subsheaf of ${\bf
R}^0p_*\U$, we can show it vanishes if the map
$$\H^0(S,\U_t)\stackrel{\mathrm{ev}}{\longrightarrow}\H^0(Z,\U_t|_Z)$$
is generically injective. Since both spaces have dimension $gn$, we
can instead prove generic surjectivity.

\begin{lem}
For generic $Z\in\mathrm{Hilb}^gS$ and $t\in S^{\prime}$ the
evaluation map
$$\H^0(S,\U_t)\stackrel{\mathrm{ev}}{\longrightarrow}\H^0(Z,\U_t|_Z)$$
is surjective.
\end{lem}

\begin{prf}
Begin with the $gn$ sections of $\U_t$. Pick a point $z_1\in S$
such that the $gn$ sections generate the fibre $(\U_t)_{z_1}$, and
remove any $n$ sections which generate $(\U_t)_{z_1}$. Now pick a
second point $z_2\in S$ such that the remaining $(g-1)n$ sections
generate the fibre $(\U_t)_{z_2}$, and remove any $n$ sections which
generate $(\U_t)_{z_2}$. Continuing in this manner will lead to
$Z:=\{z_1,\ldots,z_g\}$ such that
$$\H^0(S,\U_t)\stackrel{\mathrm{ev}}{\longrightarrow}\H^0(Z,\U_t|_Z)$$
is surjective, unless at some stage we cannot find a point $z\in S$
such that $(\U_t)_z$ is generated by $n$ independent sections of
$\U_t$. 

This would imply that these $n$ sections instead generate a
subsheaf $\cal E$ of rank $r<n$. Choose $r$ of the $n$ sections which
generate the generic fibre of $\cal E$. Thus there is an injection
$$\O^{\oplus r}\hookrightarrow{\cal E}\rightarrow{\cal
  F}$$
whose cokernel $\cal F$ is supported on a curve in $S$. Thus
$$c_1({\cal E})=c_1({\cal F})>0$$
if $\cal F$ is non-trivial. However, by stability the slope of $\cal
E$ must be less than the slope $\frac{1}{n}C^2$ of $\U_t$, and hence
can be at most zero. This implies that $c_1({\cal E})=0$, $\cal F$
vanishes, and ${\cal E}\cong\O^{\oplus r}$, which contradicts the fact
that $\cal E$ has at least $n$ sections.

Since surjectivity of $\mathrm{ev}$ is a Zariski open condition, we
conclude that $\mathrm{ev}$ is surjective for generic $t$ and $Z$.
\end{prf}



We have shown that for generic $Z\in\mathrm{Hilb}^gS$ the cohomology
of the Fourier-Mukai transform $\Psi({\cal I}_Z)$ of ${\cal I}_Z$ is
the twisted sheaf
$${\bf R}^1p_*(\U\otimes q^*{\cal I}_Z)$$
in degree one. Therefore $\Psi$ composed with taking cohomology and
the shift functor induces a rational map
$$\mathrm{Hilb}^gS=M_{(1,0,1-g)}(S)\dashrightarrow
M_{(0,D,-k)}(S^{\prime},\alpha).$$
Since both $\mathrm{Hilb}^gS$ and $M_{(0,D,-k)}(S^{\prime},\alpha)$
are smooth compact holomorphic symplectic manifolds of the same
dimension $2g$, they must be birational. In the next section we will
prove that this map is actually an isomorphism.


\section{An isomorphism of moduli spaces}

To prove that $\mathrm{Hilb}^gS$ and $M_{(0,D,-k)}(S^{\prime},\alpha)$
are isomorphic we will apply the Fourier-Mukai transform in the
reverse direction. Define $\Phi$ to be the Fourier-Mukai transform
$$\Phi^{{\cal U}^{\vee}}_{S^{\prime}\rightarrow S}:\D
(S^{\prime},\alpha)\rightarrow \D(S).$$
Note that this differs from $\Phi$ in Section 3 as we've replaced
$\cal U$ with ${\cal U}^{\vee}$, so that this new $\Phi$ is the
quasi-inverse of $\Psi$ from Section 4. We will show that the
composition of $\Phi$ with taking cohomology and the shift functor
induces an isomorphism
$$M_{(0,D,-k)}(S^{\prime},\alpha)\stackrel{\cong}{\longrightarrow}\mathrm{Hilb}^gS.$$ 
This map is clearly the inverse of the rational map from the
previous section (which is therefore also an isomorphism).

The argument we present in this section is due to Kota Yoshioka; he
kindly allowed the author to include it here.

\subsection{Preliminaries}

Since we intend to apply the Fourier-Mukai transform in the opposite
direction, we need to understand $S$ as a moduli space of twisted
sheaves on $S^{\prime}$. The $p^*\alpha-$twisted universal sheaf $\U$
on $S\times S^{\prime}$ is locally free of rank $n$. Recall that we
defined $\U_s$ to be the restriction $\U|_{s\times S^{\prime}}$ for
$s\in S$. Choosing $G$ to be $\U_s$ for some $s\in S$, as in Section
3, the degree of $\U_s$ is
$$\mathrm{deg}_G(\U_s)=0$$
where we have used the polarization $nD$ of $S^{\prime}$. Using this
definition of degree of twisted sheaves, the definition of
$\mu-$stability is the same as in the untwisted case.

\begin{lem}
The $\alpha-$twisted sheaf $\U_s$ on $S^{\prime}$ is $\mu-$stable.
\end{lem}

\begin{prf}
That $\U_s$ is stable follows from Corollary 4.5 of
Yoshioka~\cite{yoshioka04}. Furthermore, since $\alpha$ is of order
$n$ there cannot exist any locally free $\alpha-$twisted sheaves of
rank less than $n$. In particular, there cannot exist any
destabilizing sheaves of the rank $n$ twisted sheaf $\U_s$, which is
therefore $\mu-$stable.
\end{prf}

We denoted the twisted Mukai vector of $\U_s$ by $(n,E,l)$. According
to the lemma, $S$ is a moduli space of $\mu-$stable $\alpha-$twisted
sheaves on $S^{\prime}$ with Mukai vector $(n,E,l)$. Note that the
cohomological Fourier-Mukai transform $\Phi_*$ must take $(n,E,l)$ to
$(0,0,1)$.

\subsection{The proof of Theorem~\ref{main_thm}}

Let $\cal F$ be an $\alpha-$twisted sheaf on $S^{\prime}$ with twisted
Mukai vector $(0,D,-k)$. The cohomological twisted Fourier-Mukai
transform $\Phi_*$ takes $(0,D,-k)\in\tilde{H}(S^{\prime},B,\Z)$ to
$-v=-(1,0,1-g)\in\tilde{H}(S,\Z)$. Thus $\Phi$ takes $\cal F$ to an
object in $\D(S)$ with Mukai vector $-v$. As before, the image is a
priori a complex of sheaves, and the cohomology of this complex is
${\bf R}^iq_*(\U^{\vee}\otimes p^*{\cal F})$, which vanishes for $i<0$
and $i>2$. Moreover, since $\cal F$ is supported on a curve in
$S^{\prime}$, we have
$$\H^2(S^{\prime},\U_s^{\vee}\otimes{\cal F})=0$$
for all $s\in S$, and therefore ${\bf R}^2q_*(\U^{\vee}\otimes p^*{\cal
  F})$ also vanishes.

Next we will show that
$$\H^0(S^{\prime},\U_s^{\vee}\otimes{\cal F})=\mathrm{Hom}(\U_s,{\cal
  F})$$
vanishes except at a finite number of points $s\in S$. Note that since
$D$ generates $\mathrm{Pic}S^{\prime}$ over $\Z$, every curve in the
linear system $|D|$ is reduced and irreducible. This implies that
$\cal F$, whose support is a curve in $|D|$, must be $\mu-$stable as
there cannot exist a destabilizing subsheaf. The degree of $\cal F$ is 
$$\mathrm{deg}_G({\cal F})=\mathrm{deg}(\U_s^{\vee}\otimes{\cal F}).$$
Since the twisted Mukai vector of $\U_s$ is $(n,E,l)$ and that of
$\cal F$ is $(0,D,-k)$, the tensor product $\U_s^{\vee}\otimes{\cal
  F}$ will have degree $n^2D^2$, with respect to the polarization $nD$
of $S^{\prime}$.

\begin{lem}
The group $\mathrm{Hom}(\U_s,{\cal F})$ vanishes for all but finitely
many points $s\in S$. 
\end{lem}

\begin{prf}
The proof is essentially the same as the proof of part (1) of
Proposition 3.5 in Yoshioka~\cite{yoshioka01}. If
$\mathrm{Hom}(\U_s,{\cal F})$ is non-trivial then the evaluation map
$$\mathrm{Hom}(\U_s,{\cal F})\otimes \U_s\rightarrow {\cal F}$$
is surjective in codimension one. Suppose that
$\mathrm{Hom}(\U_s,{\cal F})$ is non-trivial for $N$ points
$s=s_1,\ldots,s_N$ and write $\U_i$ for $\U_{s_i}$. We claim that
$\mathrm{ker}\phi$ is a $\mu-$stable $\alpha-$twisted sheaf, where
$$\phi:\bigoplus_i \mathrm{Hom}(\U_i,{\cal F})\otimes U_i
\rightarrow{\cal F}$$
is the evaluation map.

Firstly, $\mathrm{ker}\phi$ has rank $Nn$ and 
$$\mathrm{deg}_G(\mathrm{ker}\phi)=-\mathrm{deg}_G({\cal
  F})=-n^2D^2.$$
Suppose that $\cal H$ is an $\alpha$-twisted locally free sheaf which
destabilizes $\mathrm{ker}\phi$. The rank of $\cal H$ is less than 
$Nn$ and must be a multiple of $n$, the order of $\alpha$. Together
with the fact that $\mathrm{Pic}S^{\prime}\cong\Z D$, this implies
that
$$\mathrm{deg}_G({\cal H})=\mathrm{deg}(\U_s^{\vee}\otimes{\cal H})$$
is a multiple of $n^2D^2$ (compare with the proof of
Lemma~\ref{Picard}, which also shows that every $\alpha-$twisted sheaf
on $S^{\prime}$ must have degree a multiple of $C^2=n^2D^2$). Since
$$\mu(\mathrm{ker}\phi)=\frac{-n^2D^2}{nN}<\mu({\cal
  H})=\frac{\mathrm{deg}_G({\cal H})}{\mathrm{rank}{\cal H}}$$
we conclude that $\mathrm{deg}_G({\cal H})\geq 0$.

Thus we have
$${\cal H}\hookrightarrow \bigoplus_i\mathrm{Hom}(\U_i,{\cal
  F})\otimes\U_i$$
with $\mu({\cal H})\geq 0$. Since the $\U_i$ are $\mu-$stable with
slope zero, $\cal H$ must be a direct sum
$$\U_{i_1}\oplus\U_{i_2}\oplus\cdots\oplus\U_{i_m}$$
for some $1\leq i_1\leq\ldots\leq i_m\leq N$, but this contradicts the
fact that $\mathrm{Hom}(\U_i,{\cal F})$ is non-trivial for $1\leq
i\leq N$. This establishes the claim that $\mathrm{ker}\phi$ is
$\mu-$stable.

We have an exact sequence
$$0\rightarrow\mathrm{ker}\phi\rightarrow\bigoplus_i\mathrm{Hom}(\U_i,{\cal
  F})\otimes \U_i\rightarrow {\cal F}\rightarrow\O_Z\rightarrow 0$$
where $\O_Z$ is the structure sheaf of some zero-dimensional subscheme
of the support of $\cal F$. Therefore the twisted Mukai vector of
$\mathrm{ker}\phi$ is
\begin{eqnarray*}
v^B(\mathrm{ker}\phi) & = & Nv^B(\U_s)-v^B({\cal F})+v^B(\O_Z) \\
 & = & N(n,E,l)-(0,D,-k)+(0,0,\ell(Z)) \\
\end{eqnarray*}
where $\ell(Z)$ is the length of $Z$. Thus
\begin{eqnarray*}
v^B(\mathrm{ker}\phi)^2 & = & -2N\langle (n,E,l)^{\vee},(0,D,-k)\rangle
+2N\langle (n,E,l)^{\vee},(0,0,\ell(Z))\rangle \\
 &  & +\langle (0,D,-k)^{\vee},(0,D,-k)\rangle \\
 & = & -2N-2Nn\ell(Z)+2(g-1) \\
\end{eqnarray*}
where we have used the fact that
$$\langle
(n,E,l)^{\vee},(0,D,-k)\rangle=\langle(0,0,1)^{\vee},-(1,0,1-g)\rangle=1$$
since $\Phi_*$ is an isometry. Since $\mathrm{ker}\phi$ is
$\mu-$stable, its Mukai vector must satisfy the twisted analogue of
the Bogomolov inequality, namely $v^B(\mathrm{ker}\phi)^2\geq -2$ (see
Proposition 3.6 in Yoshioka~\cite{yoshioka04}). This implies that
$N\leq g$, which completes the proof.
\end{prf}

\begin{rmk}
Suppose that $\cal F$ is supported on a curve $\iota:D\hookrightarrow
S^{\prime}$. The sheaf $\cal F$ is (generically) the push-forward
$\iota_*L$ of a line bundle $L$ on $D$. Then $S$ parametrizes a family
of rank $n$ vector bundles $\U^{\vee}_s|_D\otimes L$ on $D$. Moreover,
one can show that $\U^{\vee}_s|_D\otimes L$ has degree
$d=(g-1)n+1$. Thus there is a map from $S$ into the moduli space
$M(n,d)$ of rank $n$ degree $d$ vector bundles on $D$. In this case,
the expected codimension of the Brill-Noether locus $W^0_{n,d}\subset
M(n,d)$ is two, where a bundle $E\in W^0_{n,d}$ if $E$ has at least one
section (see Teixidor i Bigas~\cite{teixidor91}). Our lemma says that
$S\subset M(n,d)$ intersects $W^0_{n,d}$ transversely and confirms that
the Brill-Noether locus has the expected codimension. We expect that
methods like this may yield new results in higher-rank Brill-Noether
theory; one should compare our approach to Lazarsfeld's in rank
one~\cite{lazarsfeld86}.
\end{rmk}

\begin{lem}
The sheaf ${\bf R}^0q_*({\cal U}^{\vee}\otimes{\cal F})$
vanishes. 
\end{lem}

\begin{prf}
The proof is essentially the same as the proof of part (1) of
Proposition 3.5 in Yoshioka~\cite{yoshioka01}. Recall that $\cal F$ is
supported on a curve in $|D|$, which must be reduced and irreducible
since $D$ generates $\mathrm{Pic}S^{\prime}$. Moreover, $\cal F$ is of
pure dimension, i.e.\ is a torsion free sheaf on the integral
subscheme $\mathrm{Supp}{\cal F}$ of $S^{\prime}$. This implies that
${\bf R}^0q_*({\cal U}^{\vee}\otimes{\cal F})$ must be torsion free.

On the other hand, the sheaves $\cal U$ and ${\cal O}_S\boxtimes{\cal
  F}$ on $S\times S^{\prime}$ are flat over $S$, and therefore so is
${\cal U}^{\vee}\otimes{\cal F}$. Thus the previous lemma and the Base
Change Theorem imply that ${\bf R}^0q_*({\cal U}^{\vee}\otimes{\cal
  F})$ is a torsion sheaf supported in dimension zero. The lemma
follows. 
\end{prf}

Thus the Fourier-Mukai transform $\Phi(\cal F)$ of $\cal F$
has non-trivial cohomology only in degree one (i.e.\ it is WIT$^1$),
and hence $\Phi(\cal F)$ is isomorphic to the sheaf ${\bf
  R}^1q_*({\cal U}^{\vee}\otimes{\cal F})$ in degree one.

\begin{lem}
The sheaf ${\bf R}^1q_*({\cal U}^{\vee}\otimes{\cal F})$ is torsion
free.
\end{lem}

\begin{prf}
The proof is essentially the same as the proof of part (2) of
Proposition 3.5 in Yoshioka~\cite{yoshioka01}. Let $T$ be the torsion
subsheaf of ${\bf R}^1q_*({\cal U}^{\vee}\otimes{\cal F})$. By the
Base Change Theorem ${\bf R}^1q_*({\cal U}^{\vee}\otimes{\cal F})$ is
locally free over the open set 
$$\{s\in S|\H^0(S^{\prime},{\cal U}^{\vee}_s\otimes{\cal F})= 0\},$$
and thus $T$ is supported on the zero-dimensional complement of this
set.

Now apply the inverse Fourier-Mukai transform $\Psi$ to the inclusion
$$T\hookrightarrow {\bf R}^1q_*({\cal U}^{\vee}\otimes{\cal F}).$$
Since $T$ is supported in dimension zero, $\Psi T$ will simply be a
sheaf (i.e.\ $T$ is IT$_0$). On the other hand, $\Psi{\bf
  R}^1q_*({\cal U}^{\vee}\otimes{\cal F})$ will be (isomorphic to) a
sheaf in degree one since $\Psi\circ\Phi{\cal F}\cong{\cal
  F}[2]$. Since the Fourier-Mukai transform preserves inclusions,
$\Psi T$ must vanish, and therefore $T$ must vanish.
\end{prf}

We can now complete the proof of Theorem~\ref{main_thm}. The
composition of the twisted Fourier-Mukai transform $\Phi$ with taking
cohomology and the shift functor takes $\cal F$ to a torsion free
sheaf on $S$ with Mukai vector $v=(1,0,1-g)$, or in other words, to
the ideal sheaf of a length $g$ zero-dimensional subscheme of
$S$. Therefore the Fourier-Mukai transform induces an isomorphism
$$M_{(0,D,-k)}(S^{\prime},\alpha)\stackrel{\cong}{\longrightarrow}\mathrm{Hilb}^gS$$
and Theorem~\ref{main_thm} follows.

\begin{rmk}
In Section 5 of Sawon~\cite{sawon03} we asked whether the base of a
Lagrangian fibration is always a linear system of curves on a K3 or
abelian surface. This is certainly the case for the fibration in our
theorem.
\end{rmk}

The pull-back of $\O(1)$ by the projection
$$\mathrm{Hilb}^gS\rightarrow |D|\cong\P^g$$
gives a nef line bundle $L$ on $\mathrm{Hilb}^gS$ which is isotropic
and primitive with respect to the Beauville-Bogomolov quadratic
form. We must have $L=(C,-nE)$.

The twisted moduli space $M_{(0,D,-k)}(S^{\prime},\alpha)$ is a torsor
over the compactified relative Jacobian of the family of curves in
$|D|$, as are the relative compactified Picard schemes
$\overline{\mathrm{Pic}}^d$ of degree $d$. However,
$M_{(0,D,-k)}(S^{\prime},\alpha)$ is never isomorphic to a Picard
scheme: the fact that they are locally isomorphic as fibrations is
because the gerbe $\alpha$ restricted to a curve is trivializable, but
we cannot trivialize $\alpha$ globally.

\begin{lem}
The twisted moduli space $M_{(0,D,-k)}(S^{\prime},\alpha)$ is never
isomorphic to a relative compactified Picard scheme
$$\overline{\mathrm{Pic}}^d({\cal D}/|D|)\cong M_{(0,D,d+1-g)}(S^{\prime}).$$
\end{lem}

\begin{prf}
The quadratic form on the Picard lattice of
$$M_{(0,D,-k)}(S^{\prime},\alpha)\cong\mathrm{Hilb}^gS$$
looks like
$$\left(\begin{array}{cc}
2(g-1)n^2 & 0 \\
0 & -2(g-1) \\
\end{array}\right).$$
On the other hand, the Picard lattice of $M_{(0,D,d+1-g)}(S^{\prime})$
is
$$\{(a,bD,c)|\langle (0,D,d+1-g),(a,bD,c)\rangle
=-(d+1-g)a+2(g-1)b=0\}.$$
It is generated by $(0,0,1)$ and $(a_0,b_0D,0)$ where
\begin{eqnarray*}
2(g-1) & = & a_0l, \\
d+1-g & = & b_0l, \\
\end{eqnarray*}
and $l=\mathrm{g.c.d.}(2(g-1),d+1-g)$. The quadratic form is therefore
$$\left(\begin{array}{cc}
0 & -a_0 \\
-a_0 & 2(g-1)b_0^2 \\
\end{array}\right)$$
which is not equivalent to the quadratic form above.
\end{prf}


\section{Fibrations on moduli spaces of sheaves}

In this final section we will outline how the method introduced in
this paper should generalize to arbitrary moduli spaces of sheaves on
a K3 surface $S$. Let $M_v(S)$ be a moduli space of stable sheaves on
$S$ of dimension
$$2g:=\langle v,v\rangle+2>2.$$
The weight two Hodge structure of $M_v(S)$ is isomorphic to
$v^{\perp}$ (see O'Grady~\cite{ogrady97}). Suppose there exists a line
bundle with square zero with respect to the Beauville-Bogomolov
quadratic form. In other words, we have
$$w\in \H^2(M_v(S),\Z)\cong v^{\perp}\subset\tilde{\H}(S,\Z)$$
with $w^2=0$. We can regard $w=(w_0,w_2,w_4)$ as an isotropic element
of $\tilde{\H}(S,\Z)$. We can assume $w$ is primitive, and that its
first component $w_0$ is non-negative. If $w_0=0$ then we would have
$w_2^2=0$, implying $S$ is an elliptic K3 surface. It is well known
that $M_v(S)$ is a Lagrangian fibration in this case (see
Friedman~\cite{friedman98}), so we can assume $w_0$ is positive.

Define a Mukai dual K3 surface $S^{\prime}$ as the moduli space
$M_w(S)$ of stable sheaves on $S$ with Mukai vector $w$. Since $w$ is
primitive and isotropic, $S^{\prime}$ really is a K3 surface (see
Mukai~\cite{mukai87}). Now as in Lemma~\ref{HodgeH2} there is an
isometric isomorphism
$$\H^2(S^{\prime},\Z)\cong w^{\perp}/w$$
(as we saw earlier, this is the case regardless of whether
$S^{\prime}$ is a fine or non-fine moduli space). Since $v$ and $w$
are orthogonal, $v$ defines a class in $\H^2(S^{\prime},\Z)$. Up to a
sign it will be represented by a smooth curve $D$ with
$$D^2=\langle v,v\rangle =2g-2.$$
Thus $D$ will be a smooth genus $g$ curve. In particular, there is a 
(possibly twisted) Fourier-Mukai transform
$$\Psi^{\U}_{S\rightarrow S^{\prime}}:\D(S)\rightarrow
\D(S^{\prime},\alpha)$$
which at the cohomological level takes $v$ to $\pm(0,D,k)$ for some
$k\in\Z$. As in the case of the Hilbert scheme, we expect that this
Fourier-Mukai transform composed with taking cohomology and a shift
functor will induce an isomorphism of moduli spaces
$$M_v(S)\cong M_{(0,D,k)}(S^{\prime},\alpha)$$
and therefore a Lagrangian fibration
$$M_v(S)\rightarrow |D|\cong\P^g.$$

Similar arguments should apply to the generalized Kummer varieties.

\begin{flushleft}
Department of Mathematics\hfill sawon@math.sunysb.edu\\
SUNY at Stony Brook\hfill www.math.sunysb.edu/$\sim$sawon\\
Stony Brook NY 11794-3651\\
USA\\
\end{flushleft}

\end{document}